\documentclass[12pt,draft]{article}

\usepackage{amsfonts}
\usepackage{amsmath}
\usepackage{amssymb}
\usepackage{amscd}
\usepackage{amsthm}
\usepackage{indentfirst}
\usepackage[hmargin=3.6cm,vmargin=3.7cm]{geometry}

\usepackage{epsfig}

\newtheorem{thm}{Theorem}[section]
\newtheorem{lemma}[thm]{Lemma}

\newtheorem{coroll}[thm]{Corollary}

\theoremstyle{definition}

\newtheorem{defin}[thm]{Definition}
\newtheorem{rem}[thm]{Remark}

\newtheorem{question}[thm]{Question}
\newtheorem*{acknow}{Acknowledgements}
\newtheorem*{prf}{Proof}

\newcommand{\R}{{\mathbb{R}}}

\newcommand{\Z}{{\mathbb{Z}}}
\newcommand{\N}{{\mathbb{N}}}

\newcommand{\cC}{{\mathcal{C}}}

\newcommand{\fc}{{:\ }}
\newcommand{\ve}{\varepsilon}

\newcommand{\ol}{\overline}

\newcommand{\tb}{\textbf}

\DeclareMathOperator{\im}{im}

\DeclareMathOperator{\Int}{Int}

\title{On almost Poisson commutativity in dimension two}
\author{Frol Zapolsky}
\date{}

\begin{document}

\maketitle

\begin{abstract}
Consider the following question: given two functions on a symplectic manifold whose Poisson bracket is small, is it possible to approximate them in the $C^0$ norm by commuting functions? We give a positive answer in dimension two, as a particular case of a more general statement which applies to functions on a manifold with a volume form. This result is based on a lemma in the spirit of geometric measure theory. We give some immediate applications to function theory and the theory of quasi-states on surfaces with area forms.
\end{abstract}

\renewcommand{\labelenumi}{(\roman{enumi})}

\section{Introduction and results}

This note continues the theme of function theory on symplectic manifolds (albeit only in dimension two) and its relations to the theory of quasi-states, as initiated and developed, for example, in \cite{Buhovski_conv_rate_poiss_br}, \cite{Cardin_Viterbo_comm_hamiltonians}, \cite{EP_C_zero_rigidity_of_poiss_br}, \cite{EP_qs_sympl}, \cite{EPZ_qm_Poisson_br}, \cite{Zapolsky_qs_pbr_surf}, \cite{EPR_Poisson_br_qs_sympl_integr}.

Consider the following definition.

\begin{defin}
Let $M$ be a manifold of dimension $n$ and let $\Omega$ be a volume form on $M$. For $F_1,\dots,F_n \in C^\infty(M)$ define the bracket $\{F_1,\dots,F_n\} \in C^\infty(M)$ by the relation
$$\{F_1,\dots,F_n\}\Omega = dF_1\wedge \dots \wedge dF_n\,.$$
We say that the $F_i$ commute if the bracket vanishes.
\end{defin}

\begin{rem}In case $n=2$ the bracket is the same as the Poisson bracket with respect to the area form $\Omega$ (which is symplectic). Commutativity coincides with the linear dependence everywhere of the differentials $dF_i$. Although we are mainly interested in the Poisson bracket in dimension two, it makes sense to introduce this more general definition because the same method applies in order to obtain a statement which holds for the bracket on higher-dimensional manifolds as well.
\end{rem}

We use throughout the uniform, or $C^0$, norm, defined for a function $F\fc X \to \R$, where $X$ is a set, as $\|F\|:=\sup_{x \in X} |F(x)|$. For a compactly supported continuous function $F\fc M \to \R$, where $M$ is an $n$-dimensional manifold with a volume form $\Omega$ we also define the $L^1$ norm as $\|F\|_{L^1}:=\int_M |F|\Omega$. Note that for $F_1,\dots,F_n \in C^\infty(M)$ we have
$$\|\{F_1,\dots,F_n\}\|_{L^1} = \int_M |dF_1\wedge\dots\wedge dF_n|\,.$$

The main result is
\begin{thm}\label{thm_main_result}Let $M$ be a closed $n$-dimensional manifold with a volume form $\Omega$. Let $\ve \geq 0$. If $F_1,\dots,F_n \in C^\infty(M)$ satisfy
$$\|\{F_1,\dots,F_n\}\|_{L^1} \leq 2\ve\,,$$
then there are $F_1',\dots,F_n' \in C^\infty(M)$ with $\|F_i-F_i'\| \leq \ve^{1/n}$ and $\{F_1',\dots,F_n'\} \equiv 0$.
\end{thm}

\begin{rem}Note the constant $1$ before $\ve^{1/n}$. For a discussion of its sharpness see section \ref{section_disc_open_q}.
\end{rem}

Loosely rephrased, this theorem means that if $n$ smooth functions are almost commuting in the $L^1$ sense, then they can be approximated in the uniform norm by smooth functions which commute.

Let us point out some immediate consequences of this result. First, recall Cardin and Viterbo's definition of Poisson commuting continuous functions on a symplectic manifold, see \cite{Cardin_Viterbo_comm_hamiltonians}:
\begin{defin}Let $(M,\omega)$ be a symplectic manifold. Two continuous functions $F,G$ on $M$ are said to Poisson commute if there are $F_k,G_k \in C^\infty(M)$, $k\in \N$, such that $F_k \to F$, $G_k \to G$ and $\{F_k,G_k\} \to 0$ as $k \to \infty$, all in the uniform norm.
\end{defin}

We have
\begin{coroll}\label{coroll_Poisson_comm_cont_fcns}Let $(M,\omega)$ be a closed surface with an area form. Then two continuous functions $F,G \fc M \to \R$ Poisson commute if and only if there are $F_k,G_k \in C^\infty(M)$, $k \in \N$, such that $F_k\to F$, $G_k \to G$, as $k \to \infty$, in the uniform norm, and $\{F_k,G_k\} \equiv 0$ for all $k$.
\end{coroll}
That is, two continuous functions on a closed two-dimensional symplectic manifold Poisson commute if and only if they can be approximated, in the uniform norm, by Poisson commuting smooth functions.

To state the next corollary, we need to recall the notion of a quasi-state, due to Aarnes, \cite{Aarnes_quasi-states}. The reader is also referred to \cite{EPZ_qm_Poisson_br}, \cite{Zapolsky_qs_pbr_surf}, \cite{EPR_Poisson_br_qs_sympl_integr} for a connection with function theory on symplectic manifolds.
\begin{defin}
If $Z$ is a compact (Hausdorff) space, let $C(Z)$ denote the Banach algebra of all real-valued continuous functions on $Z$. Denote by $C(F)$ the closed subalgebra of $C(Z)$ generated by $F$, that is $C(F)=\{\phi\circ F\,|\,\phi\in C(\im F)\}$. A functional $\eta \fc C(Z) \to \R$ is called a quasi-state if it satisfies
\begin{enumerate}
\item $\eta(1) = 1$;
\item $\eta(F) \geq 0$ for $F \geq 0$;
\item for each $F \in C(Z)$ the restriction $\eta|_{C(F)}$ is linear.
\end{enumerate}
\end{defin}

In \cite{EP_qs_sympl}, Entov and Polterovich show that if $(M,\omega)$ is a closed surface with an area form, then a quasi-state on $M$ is linear on Poisson commutative subspaces of $C^\infty(M)$. We combine their result with corollary \ref{coroll_Poisson_comm_cont_fcns} to obtain
\begin{coroll}\label{coroll_qs_linear_comm_cont_fcns}Let $(M,\omega)$ be a closed surface with an area form. Then a quasi-state on $M$ is linear on Poisson commuting subspaces of $C(M)$.
\end{coroll}

Theorem \ref{thm_main_result} will be proved using the following lemma, which is of independent interest. First, for a map $\phi \fc \R^n \to \R^n$ define the $i$-th displacement function $\Delta_i\phi \fc \R^n \to [0,\infty)$, where $i = 1,\dots,n$, by $\Delta_i\phi(x) = |p_i(x) - p_i(\phi(x))|$, $p_i \fc \R^n \to \R$ being the projection on the $i$-th coordinate.
\begin{lemma}\label{lemma_geom_meas_thry}Let $K \subset \R^n$ be a compact set of measure $\leq \ve$. Then there is a smooth map $\phi \fc \R^n \to \R^n$ such that the displacement functions satisfy $\|\Delta_i\phi\| \leq \ve^{1/n}$, $i = 1,\dots,n$, and $\phi(K)$ has measure zero.
\end{lemma}

\begin{acknow}I would like to thank Barney Bramham for useful discussions, and Marco Mazzucchelli for listening to the preliminary version of the results, kindly proofreading the manuscript, and for useful comments.
\end{acknow}

\section{Proofs}

We begin by proving theorem \ref{thm_main_result}, assuming lemma \ref{lemma_geom_meas_thry}.
\begin{prf}[of theorem \ref{thm_main_result}]
Consider the evaluation map $\alpha \fc M \to \R^n$, $\alpha(x) = (F_1(x),\dots,F_n(x))$. Define $n_\alpha \fc \R^n \to \N \cup \{\infty\}$ by $n_\alpha(z) = \#\alpha^{-1}(z)$. The area formula (see for example \cite[theorem 3.2.3]{Federer_geom_meas}) states that $n_\alpha$ is almost everywhere real-valued and moreover
$$\int_{\R^n}n_\alpha\Omega_0 = \int_M \alpha^*\Omega_0\,.$$
Here $\Omega_0 = dx_1\dots dx_n$ is the standard density\footnote{A density on an $n$-dimensional manifold $M$ is a section of the bundle $\Lambda^nT^*M\otimes o(M)$, where $o(M)$ is the orientation line bundle of $M$.} on $\R^n$ and $\alpha^*\Omega_0$ is the pull-back density on $M$. Now $\alpha^*\Omega_0 = |dF_1\wedge \dots \wedge dF_n|$ and so
$$\int_{\R^n}n_\alpha\Omega_0 = \int_M |dF_1\wedge \dots \wedge dF_n| = \|\{F_1,\dots,F_n\}\|_{L^1}\,.$$
Denote $K = \im \alpha$. It is a compact subset of $\R^n$. Since $M$ is closed and $\R^n$ is non-compact, the degree of $\alpha$ is zero, hence zero modulo $2$, which means that $n_\alpha \geq 2$ almost everywhere on $K$. Consequently we obtain
$$2|K| = 2\int_{\im \alpha}\Omega_0 \leq \int_{\R^n} n_\alpha\Omega_0 = \|\{F_1,\dots,F_n\}\|_{L^1} \leq 2\ve\,,$$
where $|\cdot|$ is the Lebesgue measure. This shows that $|K| \leq \ve$. Lemma \ref{lemma_geom_meas_thry} yields a smooth map $\phi \fc \R^n \to \R^n$ with $\|\Delta_i\phi\| \leq \ve^{1/n}$ for all $i$ and $|\phi(K)|=0$. Define $\alpha' = \phi \circ \alpha \fc M \to \R^n$ and $F_i' = p_i \circ \alpha' \fc M \to \R$. Since $\|\Delta_i\phi\| \leq \ve^{1/n}$ for all $i$, we see that
\begin{align*}
\|F_i-F_i'\| &= \sup_M|F_i-F_i'|\\
&= \sup_M|p_i \circ \alpha - p_i \circ \phi \circ \alpha|\\
&= \sup_{\im \alpha}|p_i - p_i \circ\phi|\\
&\leq \sup_{\R^n}|p_i - p_i \circ\phi|\\
&= \|\Delta_i\phi\| \leq \ve^{1/n}\,.
\end{align*}
Moreover, since $\im \alpha' = \phi(K)$ has measure zero, the $dF_i'$ are everywhere linearly dependent, and so $\{F_1',\dots,F_n'\} \equiv 0$, as required. \qed
\end{prf}

We now prove corollary \ref{coroll_Poisson_comm_cont_fcns}.
\begin{prf}
The ``if'' part being clear, let us show the ``only if'' part. Without loss of generality assume $\int_M \omega = 1$. Suppose $F,G \in C(M)$ Poisson commute, so that there are $F_k,G_k \in C^\infty(M)$ for $k \in \N$ with $F_k \to F$, $G_k \to G$, $\{F_k,G_k\}\to 0$ in the uniform norm as $k \to \infty$. Denote $\ve_k = \frac 1 2\|\{F_k,G_k\}\|$. Then
$$\|\{F_k,G_k\}\|_{L^1} = \int_M|\{F_k,G_k\}|\omega \leq 2\ve_k\,.$$
Theorem \ref{thm_main_result} provides smooth functions $F_k',G_k'$ with $\|F_k-F_k'\|, \|G_k-G_k'\| \leq \sqrt{\ve_k}$ and $\{F_k',G_k'\} \equiv 0$. Now as $k \to \infty$,
$$\|F-F_k'\| \leq \|F-F_k\| + \|F_k-F_k'\| \leq \|F-F_k\| + \sqrt{\ve_k} \to 0\,,$$
and similarly for the $G_k'$. Thus $F_k' \to F,G_k' \to G$ as $k \to \infty$ in the uniform norm, and $\{F_k',G_k'\}\equiv 0$ for all $k$, as claimed. \qed
\end{prf}

For the proof of corollary \ref{coroll_qs_linear_comm_cont_fcns} recall that a quasi-state $\eta$ is Lipschitz with respect to the uniform norm, that is $|\eta(F) - \eta(G)| \leq \|F-G\|$ for continuous $F,G$, see \cite{Aarnes_quasi-states}.

\begin{prf}[of corollary \ref{coroll_qs_linear_comm_cont_fcns}]Denote the quasi-state by $\eta$. A quasi-state being homogeneous by definition, it suffices to show its additivity on Poisson commuting continuous functions. Thus let $F,G \in C(M)$ Poisson commute. Corollary \ref{coroll_Poisson_comm_cont_fcns} says there are $F_k,G_k \in C^\infty(M)$ such that $\{F_k,G_k\} \equiv 0$ for all $k$ and $F_k \to F, G_k \to G$ as $k \to \infty$ in the uniform norm. We have
$$|\eta(F+G) - \eta(F) - \eta(G)| = \lim_{k \to \infty}|\eta(F_k+G_k) - \eta(F_k) - \eta(G_k)| = 0\,,$$
where the first equality is due to the fact that $\eta$ is Lipschitz, while the second follows from the aforementioned result of Entov and Polterovich that a quasi-state on $M$ is linear on Poisson commuting subspaces of $C^\infty(M)$. \qed
\end{prf}

Introduce some notation. Let $\|\cdot \|$ be the Euclidean norm on $\R^n$. For $p \in \R^n$ and $\delta > 0$ let $B(p,\delta) \subset \R^n$ denote the open Euclidean ball of radius $\delta$ centered at $p$. For $\nu \in \Z^n$ we denote $C_\nu = \prod_{i=1}^n[\nu_i,\nu_i+1] \subset \R^n$, and call any such set an integer cube; also define $m_\nu = (\nu_1+\frac 1 2, \dots, \nu_n + \frac 1 2)$, which is the center of $C_\nu$.

For the proof of lemma \ref{lemma_geom_meas_thry} we need the following technical result.

\begin{lemma}\label{lemma_technical}For $\ve \in (0,\frac 1 6]$ there is a smooth map $\psi \fc \R^n \to \R^n$ which sends every integer cube to itself, is the identity on $\bigcup_{\nu \in \Z^n}B(m_\nu,\ve)$, and for every $\nu \in \Z^n$ maps $C_\nu - B(m_\nu, 2\ve)$ onto $\partial C_\nu$.
\end{lemma}

\begin{prf}[of lemma \ref{lemma_geom_meas_thry} assuming lemma \ref{lemma_technical}]
If $K$ has measure zero, the identity map does the job. Otherwise let $\gamma = |K|^{-1/n}$ and let $m_\gamma \fc \R^n \to \R^n$ be the dilation by $\gamma$, $m_\gamma(x) = \gamma x$. We have $|m_\gamma(K)| = 1$. Suppose we proved the claim of the lemma for sets of measure $1$, and let $\phi'$ be a map corresponding to $m_\gamma(K)$. Then $\phi=m_\gamma^{-1}\phi'm_\gamma$ satisfies the requirements of the lemma for $K$. Hence there is no loss of generality in assuming $|K|=1$.

If $K = C_\nu$ is an integer cube, we let $\Phi$ denote the time-$1$ map of the flow of the smooth vector field $X$ defined by $X(x) = 2\sqrt n \sigma(\|x-\nu\|)(\nu-x)/\|x-\nu\|$, where $\sigma \fc [0,\infty) \to [0,1]$ is a smooth function such that $\sigma|_{[0,\frac 1 {10}]} =0$, $\sigma|_{[\sqrt n + 2,\infty]} =0$, $\sigma|_{[\frac 1 9, \sqrt n + 1]}=1$. Then $\Psi(K) \subset C_\nu$ and $\Psi(K)$ avoids $\bigcup_{\nu' \in \Z^n}B(m_{\nu'},\frac 1 3)$.

Otherwise let $\Phi$ be the smooth map defined as follows. Denote by $\cC$ the collection of integer cubes meeting $K$, and for $C \in \cC$ let $\nu_C \in \Z^n$ be the unique integer $n$-tuple such that $C = C_{\nu_C}$. Since $K$ has measure $1$ and is not an integer cube, for each $C\in\cC$ there is $p_C \in \Int C$ and $\ve_C > 0$ such that $\ol{B(p_C,2\ve_C)} \subset \Int C$ and $B(p_C,2\ve_C) \cap K = \varnothing$. Let $\ve = \min(\frac 1 6, \min_{C \in \cC}\ve_C)$. Let $Z_C = \bigcup_{t\in[0,1]}B(tp_C+(1-t)m_{\nu_C},2\ve)$. Define the constant vector field $X_C$ on $Z_C$ via $X_C = m_{\nu_C} - p_C$ and extend it to a smooth field, still denoted by $X_C$, on $\R^n$ with compact support in $\Int C$. Let $X = \sum_{C \in \cC}X_C$ and let $\Phi$ be the time-$1$ map of the flow of $X$. Then $\Phi$ maps $B(p_C,2\ve)$ isometrically onto $B(m_{\nu_C}, 2\ve)$, and $\Phi(K) \cap \bigcup_{\nu \in \Z^n}B(m_\nu,2\ve) = \varnothing$.

Let $\psi$ be a map guaranteed by lemma \ref{lemma_technical} for $\ve$ defined as above, and put $\phi = \psi \circ \Phi$. It is easy to see that $\phi$ satisfies the requirements of the lemma.\qed
\end{prf}

Now it only remains to prove lemma \ref{lemma_technical}.

\begin{prf}
The required map is constructed in stages.

Start with a smooth map $a \fc [0,1]\to[0,1]$ which coincides with the identity map near $\frac 1 2$ and whose derivatives all vanish at $0$ and $1$. For example, $a$ can be defined by $a(t) = (t\rho(t)-1)\rho(1-t) + 1$ for $t \in (0,1)$, $a(0) = 0$, $a(1)=1$, where
$$\rho(t) = \frac{e^{-\Lambda/t}}{e^{-\Lambda/t}+e^{-\Lambda/(1-t)}}\,,$$
$\Lambda > 0$ being a sufficiently large number.

Next, define $b_n \fc [0,1]^n \to [0,1]^n$ by $b(x_1,\dots,x_n)=(a(x_1),\dots,a(x_n))$.

Now let $c \fc \partial [-\frac 1 2, \frac 1 2]^n \to \partial [-\frac 1 2, \frac 1 2]^n$ be defined as follows. If $F \subset \partial [-\frac 1 2, \frac 1 2]^n$ is an $(n-1)$-dimensional face, let $i \fc F \to [0,1]^{n-1}$ be an isometry, and let $c|_F:=i^{-1}\circ b_{n-1} \circ i$.

Let $p \fc \R^n - \{0\} \to S^{n-1}$ be the radial projection. Define $f \fc S^{n-1} \to \partial [-\frac 1 2, \frac 1 2]^n$ by $f = c \circ \big(p|_{\partial [-\frac 1 2, \frac 1 2]^n}\big)^{-1}$. Then $f$ is a smooth one-to-one and onto map from $S^{n-1}$ to $\partial [-\frac 1 2, \frac 1 2]^n$. It is a diffeomorphism when restricted to the preimage of any $(n-1)$-dimensional open face of $\partial[-\frac 1 2, \frac 1 2]^n$, and its critical values fill the complement of the union of the open faces. 

Let us construct $\psi$. Let $\lambda \fc \R \to [0,1]$ be a smooth function such that $\lambda(t) = 0$ for $t \leq \ve$, $\lambda(t) = 1$ for $t \geq 2\ve$. For $x \in C_\nu-m_\nu$ put
$$\psi(x) := \big(1-\lambda\big(\|x-m_\nu\|\big)\big)x+\lambda\big(\|x-m_\nu\|\big)\Big(m_\nu + f\Big(\frac{x-m_\nu}{\|x-m_\nu\|}\Big)\Big)\,,$$
and $\psi(m_\nu):=m_\nu$. It is an exercise to check that $\psi$ is a well-defined smooth map. Since a cube is convex, $\psi$ maps every integer cube to itself. It also follows from its definition that it is the identity on $\bigcup_{\nu \in \Z^n}B(m_\nu,\ve)$ and maps the complement of $B(m_\nu,2\ve)$ in $C_\nu$ onto $\partial C_\nu$, as required. \qed
\end{prf}

\section{Discussion and open questions}\label{section_disc_open_q}

The result stated in theorem \ref{thm_main_result} can be viewed as complementary to the so-called rigidity of Poisson brackets as shown in \cite{Buhovski_conv_rate_poiss_br}, \cite{EP_C_zero_rigidity_of_poiss_br}. Rigidity means that the functional $C^\infty \times C^\infty \to [0,\infty)$, $(F,G) \mapsto \|\{F,G\}\|$ is lower semi-continuous in the $C^0$ topology, or more informally, that it is impossible to significantly reduce the $C^0$ norm of the Poisson bracket of two smooth functions by an arbitrarily small $C^0$ perturbation. Theorem \ref{thm_main_result} means that if two functions have small Poisson bracket, the two functions can be perturbed in the $C^0$ norm so that the new functions have vanishing bracket. In view of this it is natural to ask
\begin{question}Is an analog of theorem \ref{thm_main_result} true on higher-dimensional symplectic manifolds? More precisely, given a closed symplectic manifold $(M,\omega)$ is there a constant $C>0$ such that for functions $F,G \in C^\infty(M)$ with $\|\{F,G\}\| = 1$ there are functions $F',G' \in C^\infty(M)$ such that $\|F-F'\|,\|G-G'\| \leq C$ and $\{F',G'\}\equiv 0$? If not, what kind of obstruction prevents this from happening?
\end{question}

The constant $1$ appearing as the factor before $\ve^{1/n}$ in theorem \ref{thm_main_result} is conjecturally not sharp.
\begin{question}What is the sharp constant in theorem \ref{thm_main_result}?
\end{question}
We believe that it is $\frac 1 2$. It cannot be less than $\frac 1 2$ because, as it is fairly easy to show, for any closed connected manifold $M$ of dimension $n$ there is a map $\alpha \fc M \to \R^n$ having its image equal to $[0,1]^n$ and with function $n_\alpha$ equal almost everywhere to $2$ on $\im \alpha$. The intermediate value theorem implies that any continuous map $\alpha' \fc M \to \R^n$ satisfying $\|\Delta_i(\alpha - \alpha')\| < \frac 1 2$ for all $i$ has image of positive measure. In terms of the bracket it means that there are $n$ smooth functions on $M$ with the $L^1$-norm of the bracket equal to $2 = 2\cdot 1$ such that if they are perturbed in the uniform norm by less than $\frac 1 2$, the bracket of the new functions is not identically zero.

Lemma \ref{lemma_geom_meas_thry} reminds of the classical isoperimetric inequality, in that it relates a volume measurement, that is the measure of the set, to a linear measurement, that is the maximal displacement of a smooth map contracting it to a set of measure zero. If for a compact $K \subset \R^n$ we denote $\text{thickness}\,(K) = \inf \{\max_{i=1,\dots,n}\|\Delta_i\phi\|\,|\,\phi \fc \R^n\to\R^n \text{ smooth with } |\phi(K)|=0\}$, then lemma \ref{lemma_geom_meas_thry} states that
$$\big(\text{thickness}\,(K)\big)^n \leq |K|\,.$$

Frol Zapolsky

Max Planck Institute for Mathematics in the Sciences

Inselstrasse 22

04103 Leipzig

Germany

Email: \texttt{zapolsky@mis.mpg.de}

\end{document}